\documentclass[11pt]{article}

\usepackage{amssymb,amsthm,amscd,array}
\usepackage{graphics}
\usepackage[final]{epsfig}

\let\ssection=\section
\renewcommand{\section}{\setcounter{equation}{0}\ssection}

\newcommand{\bbR}{\mathbb{R}}

\newcommand{\bbZ}{\mathbb{Z}}

\newcommand{\bbC}{\mathbb{C}}
\newcommand{\bbK}{\mathbb{K}}

\newcommand{\ad}{\mathrm{ad}}

\newcommand{\cB}{{\cal{B}}}

\newcommand{\cM}{{\mathcal{M}}}

\newcommand{\cK}{{\mathcal{K}}}
\newcommand{\cAK}{{\mathcal{AK}}}

\newcommand{\cI}{{\mathcal{I}}}
\newcommand{\mod}{\mathrm{mod}}
\newcommand{\Hom}{\mathrm{Hom}}
\newcommand{\End}{\mathrm{End}}
\newcommand{\Der}{\mathrm{Der}}

\newcommand{\sgn}{\mathrm{sign}}

\newcommand{\Alt}{\mathrm{Alt}}
\newcommand{\osp}{\mathrm{osp}}

\newcommand{\half}{\frac{1}{2}}
\newcommand{\fg}{\mathfrak{g}}

\newcommand{\fa}{\mathfrak{a}}

\chardef\s=110
\chardef\g=103

\begin{document}
\newtheorem{theorem}{Theorem}
\newtheorem{lemma}{Lemma}[section]
\newtheorem{cor}[lemma]{Corollary}
\newtheorem{conj}[lemma]{Conjecture}
\newtheorem{proposition}[lemma]{Proposition}
\newtheorem{rmk}[lemma]{Remark}
\newtheorem{exe}[lemma]{Example}
\newtheorem{defi}[lemma]{Definition}

\def\a{\alpha}
\def\b{\beta}
\def\d{\delta}
\def\e{\varepsilon}
\def\i{\iota}
\def\g{\gamma}
\def\L{\Lambda}
\def\om{\omega}
\def\Om{\Omega}
\def\r{\rho}
\def\s{\sigma}
\def\t{\tau}
\def\vfi{\varphi}
\def\vr{\varrho}
\def\l{\lambda}
\def\m{\mu}

\title{Alternated Hochschild Cohomology}

\author{Pierre Lecomte
\and
Valentin Ovsienko}

\date{}

\maketitle

\begin{abstract}
In this paper we construct a graded Lie algebra on the space of cochains on
a $\bbZ_2$-graded vector space that are skew-symmetric in the odd variables.
The Lie bracket is obtained from the classical Gerstenhaber bracket
by (partial) skew-symmetrization; the coboundary operator is  a skew-symmetrized
version of the Hochschild differential.
We show that an order-one element $m$ satisfying the zero-square condition $[m,m]=0$
defines an algebraic structure called ``Lie antialgebra'' in \cite{Ovs}.
The cohomology (and deformation) theory of these algebras is then defined.
We present two examples of non-trivial cohomology classes
which are similar to the celebrated Gelfand-Fuchs and Godbillon-Vey classes.
\end{abstract}

\maketitle


\medskip

\thispagestyle{empty}


\section{Introduction}

Let $V=V_0\oplus{}V_1$ be a $\bbZ_2$-graded vector space.
We consider the space of parity preserving
multilinear maps
\begin{equation}
\label{TwoBilM}
\vfi:\left(V_0\times\cdots\times{}V_0\right)
\otimes\left(V_1\times\cdots\times{}V_1 \right)
\to{}V,
\end{equation}
that are skew-symmetric on the subspace $V_1$.
We will define a natural structure of graded Lie algebra on this space
and develop a cohomology theory.

The graded Lie algebra on the space of all multilinear maps
on a (multi-graded) vector space $V$ is the classic Gerstenhaber algebra \cite{Ger,Ger1}.
The graded Lie algebra on the space of skew-symmetric maps
$V\wedge{}\cdots\wedge{}V\to{}V$
 is another classic graded Lie algebra called the
 Nijenhuis-Richardson algebra  \cite{NR, NR1}.
A natural homomorphism
between the Gerstenhaber algebra and that of Nijenhuis-Richardson
is given by the skew-symmetrization, see \cite{LMS}.

In this article, we introduce a graded Lie algebra  
defined by the Gerstenhaber bracket skew-symmetrized
only in a part of variables.
In this sense, our graded Lie algebra is a kind of ``intermediate
form'' between the Gerstenhaber and Nijenhuis-Richardson algebras.

The related cohomology is defined in a usual way.
We consider a \textit{parity preserving} bilinear map $m:V\times{}V\to{}V$
of the form (\ref{TwoBilM}), understood as an \textit{odd} element (of parity 1) in our graded
Lie algebra.
We assume that $m$ satisfies the condition
\begin{equation}
\label{SQRTCond}
\left[
m,\,m
\right]=0.
\end{equation}
This defines a coboundary operator
\begin{equation}
\label{DeltAd}
\d=\ad_m
\end{equation}
and the corresponding cochain complex.
Note that the condition $\d^2=0$ is an immediate consequence of
(\ref{SQRTCond}) and the graded Jacobi identity.
We then calculate the explicit combinatorial formula of the differential,
it turns out to be an amazing ``interpolation'' between
the Hochschild and Chevalley-Eilenberg differentials.

Let us stress the fact that
most of the known algebraic structures, such as associative or Lie algebras,
Lie bialgebras, Poisson structures, etc.
can be represented in terms of an order-one element $m$ of a graded Lie
algebra (usually the algebra of derivations of an associative
algebra of tensors) that satisfies the condition (\ref{SQRTCond}).
This general idea goes back to Gerstenhaber
and Nijenhuis-Richardson and became a powerful tool
to produce new (or better understand known) algebraic structures,
see \cite{LR} as an example of such approach.

It turns out quite remarkably, that a bilinear map $m$,
symmetric on $V_0$ and skew-symmetric on $V_1$,
satisfying the condition (\ref{SQRTCond}) is precisely the structure of
\textit{Lie antialgebra} introduced in~\cite{Ovs} and further studied in \cite{LM}.
This class of algebras is a particular class of Jordan superalgebras
closely related to the Kaplansky superalgebras defined in \cite{McC}.
Lie antialgebras appeared in symplectic geometry, see \cite{Ovs}.
Deducing this algebraic structure directly from the Gerstenhaber algebra 
explains its cohomologic nature.

In this paper, we define cohomology of Lie antialgebras.
We pay a special attention to lower degree cohomology spaces
and explain their algebraic sense.
In particular we show that the second cohomology space classifies extensions of Lie antialgebras
already considered in \cite{Ovs},
while the first cohomology space classifies extensions of modules.
In the end of the paper, we present two examples of non-trivial cohomology classes
generalizing two celebrated cohomology classes
of infinite-dimensional Lie algebras.
One of them is analog of the Gelfand-Fuchs class and
the second one is analog of the Godbillon-Vey class.

\section{The graded Lie algebra}\label{FirS}

We start with a brief description of the most classical
Gerstenhaber algebra and its relation to the Nijenhuis-Richardson algebra.
We discuss in some details the case of $\bbZ_2$-graded vector space.
The results of this section are well-known, see \cite{Ger,Ger1,NR,NR1,SS,LMS},
we therefore omit the proofs.

\subsection{The classical Gerstenhaber algebra}\label{GerAl}

Given a vector space $V$,
consider the space $M(V)$ of all multilinear maps $\vfi:V\times\cdots\times{}V\to{}V$.
The standard $\bbZ_{\geq-1}$-grading on $M(V)$ is given by
$$
M(V)=\bigoplus_{k\geq-1}M^k(V)=\bigoplus_{k\geq-1}\Hom_{Vect}(V^{\otimes(k+1)},V),
$$
where $Vect$ is the category of vector spaces.

The \textit{Gerstenhaber product}
of two elements $\vfi\in M^{k}(V)$ and $\vfi'\in M^{k'}(V)$ usually denoted by,
$j_{\vfi}\,\vfi'\in M^{k+k'}(V)$ is given by
\begin{equation}
\label{GerProd}
(j_{\vfi}\,\vfi')
\left(
x_0,\ldots,x_{k+k'}
\right)=
\sum_{i=0}^{k'}(-1)^{i\,k}\,
\vfi'
\left(
x_0,\ldots,\vfi(x_i,\ldots,x_{i+k}),\ldots,x_{k+k'}
\right).
\end{equation}
The classical result of Gerstenhaber \cite{Ger} states that graded bracket
\begin{equation}
\label{GerBrack}
\left[
\vfi,\vfi'
\right]=j_{\vfi}\vfi'-(-1)^{k k'}\,j_{\vfi'}\vfi
\end{equation}
equips $M(V)$ with a structure of graded Lie algebra.
In particular, it satisfies
the graded Jacobi identity
\begin{equation}
\label{GradJc}
(-1)^{k_1{k_3}}
\left[
\vfi_1,\,[\vfi_2,\vfi_3]
\right]
+(-1)^{k_1{k_2}}
\left[
\vfi_2,\,[\vfi_3,\vfi_1]
\right]
+(-1)^{k_2{k_3}}
\left[
\vfi_3,\,[\vfi_1,\vfi_2]
\right]=0,
\end{equation}
where $\vfi_i\in{}M^{k_i}(V)$ for $i=1,2,3$.

The most conceptual way to prove this statement consists in a simple observation \cite{SS} that
the Gerstenhaber algebra is nothing but
the algebra of derivations of the associative tensor algebra~$TV^*.$
Indeed, under the assumption $\dim(V)<\infty$, one has
$$
M^k(V)\cong{}(V^*)^{\otimes{}(k+1)}\otimes{}V,
$$
as a vector space\footnote{In the infinite-dimensional case, it suffices to use the inclusion
$M^k(V)\subset{}(V^*)^{\otimes{}(k+1)}\otimes{}V,$ or consider topological arguments.}.
On the other hand, a derivation
$D\in\Der\left(TV^*\right)$ is uniquely defined
(via the Leibniz rule) by its restriction to the first-order
component $V^*$ of the algebra $TV^*$, namely,
$$
D|_{V^*}:V^*\to{}TV^*.
$$
The derivation $D$ is therefore identified with an element of
$TV^*\otimes{}V$.
One thus obtains the isomorphism 
\begin{equation}
\label{GerIdenT}
\left(M(V),\,[\,,\,]\right)
\cong
\Der\left(TV^*\right)
\end{equation}
as vector spaces.
The bracket (\ref{GerBrack})
precisely the (graded) commutator in
$\Der\left(TV^*\right)$ so that the above isomorphism
is, indeed, an isomorphism of graded Lie algebras.

The following observation will be important.
The subspace
\begin{equation}
\label{SubM}
M_+(V)=\bigoplus_{k\geq0}M^k(V)
\end{equation}
is a graded Lie subalgebra of $M(V)$.

\subsection{The parity function}\label{ClassPar}

The Gerstenhaber algebra $M(V)$ is a graded Lie algebra, in particular,
it is a \textit{Lie superalgebra}: 
$$
M(V)=M(V)_0\oplus{}M(V)_1,
$$
where 
$$
M(V)_0=\bigoplus_{k\geq0}M^{2k-1}(V)
\qquad
\hbox{and}
\qquad
M(V)_1=\bigoplus_{k\geq0}M^{2k}(V)
$$
is the even and odd part of $M(V)$, respectively.
The \textit{parity}  of a homogeneous element
$\vfi\in{}M^k(V)$ (i.e., of a $(k+1)$-linear map) is defined by
\begin{equation}
\label{NaivePar}
\bar{\vfi}=k
\qquad(\mod\,2),
\end{equation}

\begin{exe}
{\rm
Vectors $v\in{}V$ viewed as elements of $M(V)$ are \textit{odd} since  $V\cong{}M^{-1}(V)$,
linear maps from $V$ to $V$ are even, bilinear maps from $V\times{}V$ to $V$ are odd, etc.
This inversed parity of elements of $V$ is important.
}
\end{exe}

\subsection{The Nijenhuis-Richardson algebra}\label{NrAlt}

The Nijenhuis-Richardson algebra, $A(V)$,
is the graded Lie algebra defined on the space of skew-symmetric multilinear maps
\begin{equation}
\label{SkeM}
\vfi:V\wedge\cdots\wedge{}V\to{}V.
\end{equation}
The space $A(V)$ is precisely the space of derivations
$\Der(\Lambda{}V^*)$, the graded Lie bracket on $A(V)$ is defined
and similar to the Gerstenhaber bracket.
This algebra is called the Nijenhuis-Richardson algebra,
it is related to the Chevalley-Eilenberg cohomology of Lie (super)algebras, see  \cite{NR}.
We omit the explicit formul{\ae}.

It was proved in \cite{LMS} that there exists a natural homomorphism of graded Lie algebras
$$
\Alt:M(V)\to{}A(V)
$$
given by skew-symmetrization.

\subsection{Gerstenhaber algebra and Hochschild cohomology} \label{AssSec}

Let us recall the most classical case of Hochschild cohomology
in the purely even case.
A bilinear map $m:V\times{}V\to{}V$ satisfies the condition (\ref{SQRTCond}) if and
only if $m$ is an associative product on~$V$.
The linear map $\d_H:=\ad_m$ from $M^k(V)$ to $M^{k+1}(V)$ is as follows
\begin{equation}
\label{HochD}
\begin{array}{rcl}
(\d_H\,\vfi)(x_0,\ldots,x_{k+1})
&=&
m(x_0,\vfi(x_1,\ldots,x_{k+1}))\\[6pt]
&-&
\displaystyle
\sum_{i=0}^k
(-1)^i\,\vfi(x_0,\ldots,m(x_i,x_{i+1}),
\ldots,x_{k+1})\\[14pt]
&+&
\displaystyle
(-1)^{k}\,m(\vfi(x_0,\ldots,x_k),x_{k+1}),
\end{array}
\end{equation}
for an arbitrary $\vfi\in{}M^{k}(V)$.
This map coincides with the classic Hochschild differential;
the corresponding cohomology is the classic Hochschild cohomology of the
associative algebra $(V,m)$, see \cite{Lod} for more details.

\subsection{Nijenhuis-Richardson algebra and cohomology of Lie superalgebras} \label{LieASec}

A skew-symmetric bilinear map $m\in{}A^1(V)$ satisfies the condition (\ref{SQRTCond}) 
if and only if $\fg=(V,m)$ is a Lie (super)algebra. 
One now uses the notation
$$
\left[
x_1,x_2
\right]:=
m(x_1,x_2),
$$
for all $x_1,x_2\in{}V$.
The map (\ref{DeltAd}) coincides
with the Chevalley-Eilenberg differential:
\begin{equation}
\label{CED}
\begin{array}{rcl}
(\d_{CE}\,\vfi)(x_0,\ldots,x_{k+1})
&=&
\displaystyle
\sum_{i=0}^{k+1}
(-1)^{i}
\left[
x_i,\vfi(x_0,\ldots,\widehat{x_i},\ldots,x_{k+1})
\right]\\[16pt]
&&+
\displaystyle
\sum_{0\leq{}i<j\leq{}k+1}
(-1)^{i+j}\,
\vfi([x_i,x_j],x_0,\ldots,\widehat{x_i},\ldots,
\widehat{x_{j}},
\ldots,x_{k+1}),
\end{array}
\end{equation}
for every $\vfi\in{}A^{k}(V)$.

\subsection{$\bbZ_2$-graded case}\label{Z2Grad}

It is easy to generalize the above definitions in the case of
$\bbZ_2$-graded vector space $V=V_0\oplus{}V_1$.
The operations (\ref{GerProd}) and (\ref{GerBrack})
are now defined via the standard sign rule.
To make this paper self-content, we give here all the necessary details, see also \cite{LMS}.
Note that the $\bbZ_2$-graded case is a particular case of the multi-graded Gerstenhaber algebra
defined in \cite{LMS}.

The space $M^k(V)$ is again a $\bbZ_2$-graded vector space.
However, the parity function is more sophisticated than (\ref{NaivePar}).
Every element $\vfi\in{}M^{k}(V)$ is a sum of two maps
$$
\vfi=\vfi_0+\vfi_1,
$$
where $\vfi_i$ is with values in $V_i$, for $i=0,1$.
Furthermore, each of the maps $\vfi_0$ and $\vfi_1$ is a sum of
homogeneous components.

\begin{defi}
\label{ParDef}
{\rm
(i)
We call an element $\vfi\in{}M^k(V)$ \textit{homogeneous of degree} $(p,q)$,
where $p+q=k+1$, if it is of the form:
$$
\vfi:V_{i_0}\times\cdots\times{}V_{i_k}\to{}V_i
$$
where $i_0+\cdots+i_k=q$.
For instance, bilinear maps from $V_0\times{}V_1$ (or from  $V_1\times{}V_0$)
to $V_i$ are homogeneous of degree $(1,1)$, whereas the maps defined on
$V_0\times{}V_0$ and $V_1\times{}V_1$ are homogeneous of degree $(2,0)$ and $(0,2)$,
respectively.
Denote by $M^{p,q}(V)$ the space of all homogeneous maps of degree $(p,q)$.

(ii)
The \textit{parity} of a homogeneous $\vfi_i\in{}M^{p,q}(V)$ is defined by
\begin{equation}
\label{Parit}
\bar{\vfi}_i:=p+i+1
\qquad(\mod\,2),
\end{equation}
where $i=0$ or $1$.
}
\end{defi}

The Gerstenhaber product on $M(V)$ is again defined 
as an alternated sum of the terms
$$
\vfi'
\left(
x_0,\ldots,\vfi(x_i,\ldots,x_{i+k}),\ldots,x_{k+k'}
\right).
$$
However, unlike (\ref{GerProd}), the sign $(-1)^{ik}$ is replaced by 
$(-1)^{(\bar{x}_0+\cdots+\bar{x}_{i-1})\bar{\vfi}}$.
The Gerstenhaber bracket is then given by
\begin{equation}
\label{GerBrackGrad}
\left[
\vfi,\vfi'
\right]=j_{\vfi}\vfi'-(-1)^{\bar{\vfi}\,\bar{\vfi'}}\,j_{\vfi'}\vfi
\end{equation} 
is again a graded Lie algebra structure.

\begin{rmk}
{\rm
The space $V$ can again be looked at as a subspace of $M(V)$,
namely $V\cong{}M^{0,0}(V)$.
With the above definition of parity (\ref{Parit}), we have:
$\bar{x}=1$, for $x\in{}V_0$ and
$\bar{y}=0,$ for $y\in{}V_1$.
This means, the parity of elements of $V$ is inversed.
}
\end{rmk}

\section{Lie antialgebras}

The notion of Lie antialgebra was introduced in \cite{Ovs}
in the context of symplectic geometry (see also \cite{GO} for the first example)
and was further studied in \cite{MG,LM,SS}.
It was shown in~\cite{LM} that Lie antialgebras is a particular case of
Jordan superalgebras.
Lie antialgebras are closely related to Lie superalgebras.
More precisely for every Lie antialgebra $\fa$, there is a Lie superalgebra $\fg_\fa$ 
called the adjoint Lie superalgebra \cite{Ovs} acting on $\fa$.

\subsection{The definition}\label{Dex}

We give two equivalent definitions, each of them has its advantages.

\begin{defi}
{\rm
A Lie antialgebra $(\fa,\cdot)$ is a commutative $\bbZ_2$-graded algebra:
$\fa=\fa_0\oplus\fa_1$ and
$\fa_i\cdot\fa_j\subset\fa_{i+j},$
so that  for all homogeneous elements $a,b\in\fa$ one has
\begin{equation}
\label{SkewP}
a\cdot{}b\,=\,(-1)^{\bar{a}\bar{b}}
\,b\cdot{}a
\end{equation}
where $\bar{a}$ is the degree of $a$,
satisfying the third-order identities:
\begin{eqnarray}
\label{AssCommT}
x_1\cdot\left(x_2\cdot{}x_3\right)
&=&
\left(x_1\cdot{}x_2\right)\cdot{}x_3,\\[10pt]
\label{CacT}
x_1\cdot
\left(x_2\cdot{}y\right)
&=&
\textstyle
\half
\left(x_1\cdot{}x_2\right)
\cdot{}y,\\[10pt]
\label{ICommT}
x\cdot\left(y_1\cdot{}y_2\right)
&=&
\left(x\cdot{}y_1\right)\cdot{}y_2
\;+\;
y_1\cdot{}
\left(x\cdot{}y_2
\right),\\[10pt]
\label{Jack}
y_1\cdot\left(y_2\cdot{}y_3\right)
&+&
y_2\cdot\left(y_3\cdot{}y_1\right)
\;+\;
y_3\cdot\left(y_1\cdot{}y_2\right)
=0,
\end{eqnarray}
for all $x_i\in\fa_0$ and $y_i\in\fa_1$.
In particular, $\fa_0$ is a commutative associative subalgebra.
}
\end{defi}

An equivalent definition is as follows.
A $\bbZ_2$-graded commutative algebra $\fa$
is a Lie antialgebra if and only if the following three conditions are satisfied.

\begin{enumerate}
\item
The subalgebra $\fa_0\subset\fa$ is associative.

\item
For all $x_1,x_2\in\fa_0$, the operators of left multiplication commute:
$x_1\cdot\left(x_2\cdot{}a\right)=x_2\cdot\left(x_1\cdot{}a\right)$.

\item
For every $y\in{}\fa_1$, the operator of right multiplication by $y$
is an \textit{odd} derivation of $\fa$, i.e., one has
\begin{equation}
\label{InvGrad}
\left(
a\cdot{}b
\right)\cdot{}y=
\left(a\cdot{}y\right)\cdot{}b
+
(-1)^{\bar{a}}\,
a\cdot{}\left(b\cdot{}y\right).
\end{equation}
\end{enumerate}

\subsection{Main examples}\label{Dex}

\begin{exe}
{\rm
Our first example of a Lie antialgebra is 3-dimensional
algebra known as \textit{tiny Kaplansky Superalgebra} and denoted by $K_3$.
This algebra has the basis $\{\e;a,b\}$,
where $\e$ is even and $a,b$ are odd,
satisfying the relations
\begin{equation}
\label{aslA}
\begin{array}{l}
\e\cdot{}\e=\e,\\[6pt]
\e\cdot{}a=\half\,a,
\quad
\e\cdot{}b=\half\,b,\\[6pt]
a\cdot{}b=\half\,\e.
\end{array}
\end{equation}
It is simple, i.e., it contains no non-trivial ideal.
}
\end{exe}

The corresponding algebra of derivations is the simple Lie superalgebra $\osp(1|2)$
of linear symplectic transformations of a $2|1$-dimensional symplectic space.
Let us mention that the algebra $K_3$ plays an important r\^ole in
the study of some exceptional Jordan superalgebras, cf.~\cite{BE}.

\begin{exe}
\label{ExMain}
{\rm
The main example of (an infinite-dimensional) Lie antialgebra is the
\textit{conformal Lie antialgebra} $\cAK(1)$.
This is a simple infinite-dimensional Lie antialgebra with the basis
$$
\textstyle
\left\{
\e_n,\;n\in\bbZ;
\quad
a_i,\;i\in\bbZ+\half
\right\},
$$
where $\e_n$ are even and $a_i$ are odd
and satisfy the following relations
\begin{equation}
\label{GhosRel}
\begin{array}{rcl}
\e_n\cdot{}\e_m&=&
\e_{n+m},\\[8pt]
\e_n\cdot{}a_i &=&
\half\,a_{n+i},\\[8pt]
a_i\cdot{}a_j &=&
\half\left(j-i\right)\e_{i+j}.
\end{array}
\end{equation}
The algebra $\cAK(1)$ is closely related
to the well-known Neveu-Schwarz conformal Lie
superalgebra,~$\cK(1)$, namely
$\cK(1)=\mathrm{Der}(\cAK(1)),$ see~\cite{Ovs}.
Note that $\cAK(1)$ contains infinitely many copies of $K_3$, among which we quote the one with basis
$\left\{\e_0;a_{-\frac 1 2},a_{\frac 1 2}\right\}.$}
\end{exe}

A very similar algebra, called  the full derivation algebra, was considered in \cite{McC}.
This algebra is also defined by the formul{\ae} (\ref{GhosRel}), but
the odd generators $a_i$ are indexed by integer $i$'s.

\begin{exe}
\label{ExMainBis}
{\rm
Our next  example is the simple infinite-dimensional Lie antialgebra $\cM^1$ which is
a ``truncated version'' of $\cAK(1)$.
The algebra $\cM^1$ is the algebra of formal series
in the elements of the basis:
$$
\textstyle
\left\{
\a_n,
\,
n\geq0,
\quad
\,a_i,
\,
i\geq-\half
\right\},
$$
subject to the relations (\ref{GhosRel}).

We understand this algebra as analog of
the Lie algebra, $W_1$, of formal vector fields on the real line.
Note that $W_1$ plays an important role in algebra and in topology (see, e.g., \cite{Fuk}).
}
\end{exe}

Further examples of finite-dimensional Lie antialgebras can be found in \cite{Ovs}.
An interesting series of simple infinite-dimensional examples
related to punctured Riemannian surfaces is constructed in \cite{SS}.

\section{Alternated Gerstenhaber algebra}\label{ParityPres}

In this section, we give our main construction 
of graded Lie algebra that will provide the cohomology theory of Lie antialgebras.

\subsection{The space of $(p,q)$-linear maps}\label{pq}

Let $V=V_0\oplus{}V_1$ be a $\bbZ_2$-graded vector space.
We will consider $(p,q)$-linear maps of the following very special form
\begin{equation}
\label{NoMix}
\vfi:\left(V_0\times\cdots\times{}V_0\right)
\otimes\left(V_1 \times\cdots \times{}V_1 \right)
\to{}V,
\end{equation}
where $p,q$ are non-negative integers.
The vector space of $(p,q)$-linear maps of the form (\ref{NoMix}) is isomorphic to
$(V_0^*)^{\otimes{}p}\otimes(V_1^*)^{\otimes{}q}\otimes{}V$.

In other words, we do not mix $V_0$ and $V_1$, for instance, we do not consider
maps of the form
 $V_1\times{}V_0\to{}V$, or $V_0\times{}V_1\times{}V_0\to{}V$, etc.

The space of $(p,q)$-linear maps of the form (\ref{NoMix})
 can be viewed as the quotient of all $(p,q)$-linear maps
by the maps vanishing on the subspace
$V_0^{\otimes{}p}\otimes{}V_1^{\otimes{}q}\subset{}V^{\otimes(p+q)}$.

\begin{proposition}
\label{QuotPro}
The space of multilinear maps (\ref{NoMix}) is isomorphic to
the quotient $M^{p,q}(V)/\cI$, where $\cI$ is the subspace
\begin{equation}
\label{Ideal}
\cI=\{\vfi\in{}M^{p,q}(V),\;\vfi|_{V_0^{\otimes{}p}\otimes{}V_1^{\otimes{}q}}=0\}.
\end{equation}
\end{proposition}

\begin{proof}
Straighforward.
\end{proof}

\begin{rmk}
{\rm
By definition of the space $\cI$, see formula (\ref{Ideal}),
for a non-zero homogeneous element $\vfi\in\cI$ of degree $(p,q)$, one has
\begin{equation}
\label{Condit}
p+q\geq1.
\end{equation}
Indeed, if $p=q=0$, then $\vfi$ is just an element of $V$ which has to be zero.
In the sequel, we will often restrict our considerations
to the subspace of maps (\ref{NoMix}) satisfying the condition (\ref{Condit}).
}
\end{rmk}

We will denote $x_i$ the elements of $V_0$
and $y_j$ the elements of $V_1$, so that we write
$$\vfi(x_1,\ldots,x_p;\,y_1,\ldots,y_q)$$ for the evaluation of $\vfi$
on some particular arguments.

\subsection{The Gerstenhaber product}\label{GerAlt}

It turns out that the classical Gerstenhaber product on the space $M_+(V)$,
see (\ref{SubM}), can be pushed forward to the quotient-space of the maps (\ref{NoMix}).
With the above notation, one has:
$$
M_+(V)=\bigoplus_{p+q\geq1}M^{p,q}(V).
$$

\begin{proposition}
\label{IdPro}
The subspace $\cI\subset{}M_+(V)$ is an ideal with respect to the Gerstenhaber product on $M_+(V)$.
\end{proposition}

\begin{proof}
Given two elements, $\vfi$ and $\vfi'$,  of $M_+(V)$,
we have to prove that $j_\vfi\vfi'\in\cI$ whenever either of the arguments,
$\vfi$ or $\vfi'$, belongs to $\cI$.
We can assume that $\vfi$ and $\vfi'$ are $(p,q)$-linear and $(p',q')$-linear, respectively.

Every map (\ref{NoMix}) is a sum of two terms:
$$
\vfi=\vfi_0+\vfi_1,
$$
where $\vfi_0$ is with values in $V_0$ and $\vfi_1$  is with values in $V_1$.
Recall also that the parity of a homogeneous $(p,q)$-linear map is defined by  (\ref{Parit}).

We have four different cases.

(a)
If $\vfi=\vfi_0$ is a $(p,0)$-linear, i.e., a map with values in $V_0$ and such that
$q=0$, then the Gerstenhaber product $j_{\vfi}\,\vfi'$ is a $(p+p'-1,q')$-linear map.
Restricting to
$V_0^{\otimes(p+p'-1)}\otimes{}V_1^{\otimes{}q'}$ one obtains:
\begin{equation}
\label{GerProdp0}
\begin{array}{l}
\displaystyle
(j_{\vfi_0}\,\vfi')
\left(
x_1,\ldots,x_{p+p'-1};\,y_1,\ldots,y_{q'}
\right)=\\[6pt]
\displaystyle
\;\;\;\;\;\;\;\;
\sum_{i=0}^{p'-1}(-1)^{i\,p}\,
\vfi'
\left(
x_1,\ldots,\vfi_0(x_{i+1},\ldots,x_{i+p}),\ldots,x_{p+p'-1};\,y_1,\ldots,y_{q'}
\right).
\end{array}
\end{equation}
The sign $(-1)^{i\,p}$ in formula (\ref{GerProdp0})
is due to the fact that the variables $x_i$ are odd and $\bar{\vfi_0}=p$, according to (\ref{Parit}).

(b)
For a $(p,q)$-linear map $\vfi_0$ (with values in $V_0$) such that $q\not=0$,
we obtain a $(p+p'-1,q+q')$-linear map $j_{\vfi}\,\vfi'$ whose restriction to
$V_0^{\otimes(p+p'-1)}\otimes{}V_1^{\otimes{}(q+q')}$ reads
\begin{equation}
\label{GerProdGrad}
\begin{array}{l}
\displaystyle
(j_{\vfi_0}\,\vfi')
\left(
x_1,\ldots,x_{p+p'-1};\,y_1,\ldots,y_{q+q'}
\right)=\\[10pt]
\displaystyle
\;\;\;\;\;\;\;\;
(-1)^{(p'-1)\,\bar{\vfi_0}}\,
\vfi'
\left(
x_1,\ldots,x_{p'-1},
\vfi_0(x_{p'},\ldots,x_{p+p'-1};\,y_1,\ldots,y_q);\,y_{q+1},\ldots,y_{q+q'}
\right),
\end{array}
\end{equation}
where the parity $\bar{\vfi}_0$ is again calculated according to (\ref{Parit}),
i.e., $\bar{\vfi}_0=p+1$.
(One does not get any summation since $\vfi_0$ can stand in the only position.)

(c)
For a $(p,q)$-linear map $\vfi_1$ (with values in $V_1$)
such that $p\not=0$, the Gerstenhaber product $j_{\vfi}\,\vfi'$ is a $(p+p',q+q'-1)$-linear map.
Restricting to
$V_0^{\otimes(p+p')}\otimes{}V_1^{\otimes{}(q+q'-1)}$ one obtains:
\begin{equation}
\label{GerProdGradBis}
\begin{array}{l}
\displaystyle
(j_{\vfi_1}\,\vfi')
\left(
x_1,\ldots,x_{p+p'};\,y_1,\ldots,y_{q+q'-1}
\right)=\\[10pt]
\displaystyle
\;\;\;\;\;\;\;\;
(-1)^{p'\,p}\,
\vfi'
\left(
x_1,\ldots,x_{p'};\,
\vfi_1(x_{p'+1},\ldots,x_{p+p'};\,y_{1},\ldots,y_{q}),
\ldots,y_{q+q'-1}
\right),
\end{array}
\end{equation}
note that $\bar{\vfi_1}=p$.

(d)
Finally, for a $(0,q)$-linear map $\vfi_1$ (with values in $V_1$ with $p=0$)
the Gerstenhaber product $j_{\vfi}\,\vfi'$ is a $(p',q+q'-1)$-linear map that restricted
to $V_0^{\otimes{}p'}\otimes{}V_1^{\otimes{}(q+q'-1)}$ reads
\begin{equation}
\label{GerProdGradTre}
\begin{array}{l}
\displaystyle
(j_{\vfi_1}\,\vfi')
\left(
x_1,\ldots,x_{p'};\,y_1,\ldots,y_{q+q'-1}
\right)=\\[6pt]
\displaystyle
\;\;\;\;\;\;\;\;
\sum_{i=0}^{q'-1}
\vfi'
\left(
x_1,\ldots,x_{p'};\,y_{1},\ldots,y_{i},
\vfi_1(y_{i+1},\ldots,y_{i+q}),
\ldots,y_{q+q'-1}
\right).
\end{array}
\end{equation}
Note that the map $\vfi_1$ is even so that there is no sign in the above product.

In all the above cases, the Gerstenhaber product $j_{\vfi}\,\vfi'$ belongs to $\cI$
provided at least one of the maps $\vfi$ and $\vfi'$ belongs to $\cI$.
\end{proof}

It follows important statement immediately from the above proposition.

\begin{cor}
The space $(V_0^*)^{\otimes{}p}\otimes(V_1^*)^{\otimes{}q}\otimes{}V$
is a graded Lie algebra, the Lie structure being obtained by the projection
of the Gerstenhaber bracket (\ref{GerBrackGrad}) to the quotient space
$M_+(V)/\cI$.
\end{cor}

\begin{rmk}
{\rm
Everywhere in this section we adopted
the condition (\ref{Condit}).
This condition is essential,
in particular, Proposition \ref{IdPro} fails without it.
To give the simplest example, consider a $(0,0)$-homogeneous map $\vfi_0$
(i.e., an element of $V_0$) and a bilinear map $\vfi'\in{}M^{(1,1)}(V)$.
Then, the Gerstenhaber product $j_{\vfi_0}\vfi'$ is a $(0,1)$-linear map
such that
$$
j_{\vfi_0}\vfi'(y)=\vfi'(\vfi_0,y)+\vfi'(y,\vfi_0).
$$
The second term does not vanish and therefore $j_{\vfi_0}\vfi'\not\in\cI$.
}
\end{rmk}

\subsection{The graded Lie algebra $\mathfrak{al}(V)$}\label{SpaC}

We are ready to introduce our main space.

\begin{defi}
{\rm
A $(p,q)$-linear map $\vfi$ is called \textit{parity preserving} if
$$
\left\{
\begin{array}{rl}
\vfi_0=0, & \hbox{if $q$ is odd},\\[4pt]
\vfi_1=0, & \hbox{if $q$ is even},
\end{array}
\right.
$$
where $\vfi_0$ is with values in $V_0$ and $\vfi_1$ is with values in $V_1$.
}
\end{defi}

Note that, in the case of parity preserving maps, the parity function can also be calculated
via the formula:
\begin{equation}
\label{TheParF}
\bar{\vfi}=p+q+1,
\end{equation}
which is in full accordance with (\ref{Parit}) as $q+i=0\,(\mod\,2)$.

\begin{exe}
{\rm
(i) 
Elements of $V_0$ are parity preserving $(0,0)$-linear maps,
while elements of $V_1$ are not parity preserving;
every element $x\in{}V_0$ is an odd $(0,0)$-linear map, i.e., $\bar{x}=1$.

(ii)
Linear parity preserving maps send $V_0$ to $V_0$ and $V_1$ to $V_1$,
all these maps are even.

(iii)
A parity preserving bilinear map can be of one of the following three forms:
\begin{equation}
\label{MapM}
m'_0:V_0\times{}V_0\to{}V_0,
\qquad
m''_0:V_1\times{}V_1\to{}V_0,
\qquad
m_1:V_0\times{}V_1\to{}V_1.
\end{equation}
In all of these three cases, the bilinear map is odd.
}
\end{exe}

\begin{defi}
{\rm
We denote by $\mathfrak{al}(V)$ the space of parity preserving multilinear maps  (\ref{NoMix})
that are \textit{skew-symmetric in $y$-variables} and satisfy the condition (\ref{Condit}).
}
\end{defi}

There is a natural projection, $\Alt$, from the space of
all parity preserving maps satisfying (\ref{Condit}) to $\mathfrak{al}(V)$.
This projection is defined via
skew-symmetrization in $y$-variables:
\begin{equation}
\label{AltyEq}
(\Alt\,\vfi)
\left(
x_1,\ldots,x_p;\,
y_1,\ldots,y_q
\right)=
\frac{1}{q!}\,
\sum_{\s\in{}S_{q}}\sgn(\s)\,\vfi
\left(
x_1,\ldots,x_p;\,
y_{\s(1)},\ldots,y_{\s(q)}
\right).
\end{equation}

Let us define a bilinear skew-symmetric operation on $\mathfrak{al}(V)$.
Given two homogeneous elements $\vfi,\vfi'\in\mathfrak{al}(V)$,
the Gerstenhaber bracket
$[\vfi,\vfi']$
is of course not necessarily in $\mathfrak{al}(V)$.
We define the following operation:
\begin{equation}
\label{OurBr}
[\vfi,\vfi']_\mathfrak{al}:=
\Alt\,[\vfi,\vfi'].
\end{equation}
Our first main result is the following.

\begin{theorem}
\label{TheDefP}
The space $\mathfrak{al}(V)$ equipped with the bracket (\ref{OurBr})
is a graded Lie algebra.
\end{theorem}

\begin{proof}
It will suffice to check that the Gerstenhaber bracket commutes with 
the map $\Alt$, that is
\begin{equation}
\label{Alt}
\Alt\,[\vfi,\vfi']=\Alt\,[\Alt\,\vfi, \Alt\,\vfi'],
\end{equation}
for all $\vfi,\vfi'\in{}M_+(V)$.
Indeed, the Jacobi identity for $[.,.]_\mathfrak{al}$ will then follow from the Jacobi identity
for the Gerstenhaber bracket.
 
Formula (\ref{Alt}) is satisfied whenever  both $\vfi$ and $\vfi'$ are 
skew-symmetric in $y$-variables since $\Alt$ is a projection and one has
$\Alt\,\vfi=\vfi$ and $\Alt\,\vfi'=\vfi'$ in this case. 
It remains to prove this formula in the case where
either $\vfi$ or $\vfi'$ belongs to the kernel of $\Alt$.
 
Consider the Gerstenhaber product (\ref{GerProdp0}-\ref{GerProdGradTre}).

In the cases a) and b), where $\vfi=\vfi_0$, the product $j_{\vfi_0}\vfi'$
belongs to $\ker(\Alt)$ if either $\vfi_0$ or $\vfi'$ belongs to $\ker(\Alt)$.

In the cases c) and d), $j_{\vfi_1}\vfi'\in\ker(\Alt)$
whenever $\vfi_1\in\ker(\Alt)$. 
It remains to show that $j_{\vfi_1}\vfi'\in\ker(\Alt)$, if $\vfi'\in\ker(\Alt)$.  
Indeed, due to the parity preserving property, the transposition
$$
\left(\vfi_1(x_{p'+1},\ldots,x_{p+p'};\,y_i,\ldots,y_{i+q}),y_j\right)
\leftrightarrow
\left(y_j,\vfi_1(x_{p'+1},\ldots,x_{p+p'};\,y_i,\ldots,y_{i+q})\right)
$$
is an odd permutation of $\{y_i,\ldots,y_{i+q},y_j\}$. 
Therefore $\Alt(j_{\vfi_1}\vfi')$ is the skew-symmetrization of 
\[
\Alt(\vfi')(x_1,\ldots,x_p;\vfi_1(x_{p'+1},\ldots,x_{p+p'};y_1,\ldots,y_q),y_{q+1},\ldots,y_{q+q'-1})
\]
with respect to the $y$'s. Hence it vanishes.
\end{proof}

Note that the parity preserving condition is essential, identity (\ref{Alt}) fails without it.

\begin{exe}
\label{ImportEx}
{\rm
Let $m\in\mathfrak{al}(V)$ be a bilinear map, we will calculate explicitly the
bracket $[m,m]_\mathfrak{al}$ of $m$ with itself (recall that $\bar{m}=1$).
Note that $m$ is of the form
$$
m=m'_0+m''_0+m_1,
$$
cf. formula (\ref{MapM}), where $m''_0$ is skew-symmetric.
One has:
$$
\begin{array}{lll}
\half\,[m,m]_\mathfrak{al}(x_1,x_2,x_3)&=&
m'_0(m'_0(x_1,x_2),x_3)-m'_0(x_1,m'_0(x_2,x_3)),\\[6pt]
\half\,[m,m]_\mathfrak{al}(x_1,x_2;\,y)&=&m_1(m'_0(x_1,x_2);\,y)-m_1(x_1,m_1(x_2;\,y)),
\end{array}
$$
this is quite obvious since there is not skew-symmetrization.
Furthermore,
$$
\textstyle
\frac{1}{4}\,[m,m]_\mathfrak{al}(x;\,y_1,y_2)=
m''_0(m_1(x;\,y_1),y_2)-m''_0(m_1(x;\,y_2),y_1)
-2\,m'_0(x;\,m''_0(y_1,y_2)),
$$
since $m''_0$ is skew-symmetric.
Finally, one has
$$
\textstyle
\frac{1}{6}\,[m,m]_\mathfrak{al}(y_1,y_2,y_3)=
m_1(m''_0(y_1,y_2);\,y_3)+m_1(m''_0(y_2,y_3);\,y_1)+m_1(m''_0(y_3,y_1);\,y_2),
$$
according to (\ref{GerProdGrad}).
}
\end{exe}

\section{Cohomology of Lie antialgebras} \label{CohSec}

In this section, we define and calculate explicitly
the cochain operator associated to an arbitrary element
$m\in{}\mathfrak{al}^1(V)$ satisfying the condition (\ref{SQRTCond}).
This is the cohomology theory we are interested in.
We start with a brief account on the classic Hochschild cohomology
defined within the context of Gerstenhaber algebra,
as well as the Chevalley-Eilenberg cohomology of Lie (super)algebras
related to the Nijenhuis-Richardson algebra.

\subsection{Combinatorial formula for the differential on $\mathfrak{al}(V)$}\label{TheDifSec}

Consider now the graded Lie algebra $\mathfrak{al}(V)$ defined in Section \ref{SpaC}.
Let $m\in{}\mathfrak{al}^1(V)$ be a parity preserving bilinear map satisfying the condition
$[m,m]_\mathfrak{al}=0$.
In this section, we calculate the combinatorial expression of the differential
$\d=\ad_m$.

Of course, the operator $\d$ increases the degree of the cochains.
We will use the notation $\d^k$ for the restriction of $\d$ to $\mathfrak{al}^{k}(V)$, viz
$$
\d^k:\mathfrak{al}^{k}(V)\to{}\mathfrak{al}^{k+1}(V),
$$
where $k\geq0$ (in accordance with the condition (\ref{Condit})).
Since $m$ preserves the parity,
it follows that the operator $\d^k$ has three terms:
\begin{equation}
\label{TheCobOpTot}
\d^k=\d^k_{1,0}+\d^k_{0,1}+\d^k_{-1,2},
\end{equation}
where $\d^k_{i,j}:\mathfrak{al}^{p,q}(V)\to{}\mathfrak{al}^{p+i,q+j}(V)$,
for $p+q=k+1$.

The following statement follows from Theorem \ref{TheDefP}.

\begin{cor}
\label{CombThm}
The operator $\d^k=\d^k_{1,0}+\d^k_{0,1}+\d^k_{-1,2}$ defined by:

\noindent
(i)
if $q=0$, then $\d_{1,0}$ is the standard Hochschild differential,
if $q>0$, then
\begin{equation}
\label{TheCobOp10}
\begin{array}{l}
(\d^k_{1,0}\vfi)(x_0,\ldots,x_{p};\,y_0,\ldots,y_{q-1})
=\\[10pt]
\qquad
\qquad
\qquad
\displaystyle
-m\left(
{x_0},\,\vfi(x_1,\ldots,x_{p};\,y_0,\ldots,y_{q-1})
\right)\\[10pt]
\qquad
\qquad
\qquad
\displaystyle
+\sum_{i=0}^{p-1}
(-1)^i\,\vfi(x_0,\ldots,x_{i-1},m(x_i,x_{i+1}),
x_{i+2},\ldots,x_{p};\,y_0,\ldots,y_{q-1})\\[14pt]
\qquad
\qquad
\qquad
\displaystyle
+
\frac{1}{q}
\sum_{j=0}^{q-1}(-1)^{p+j}\,
\vfi(x_0,\ldots,x_{p-1};\,
m(x_{p},y_j),y_0,\ldots,\widehat{y_j},\ldots,y_{q-1});
\end{array}
\end{equation}

\noindent
(ii)
if $p=0$ and $q$ is odd, then $\d^k_{0,1}$ is given by
\begin{equation}
\label{TheCobOp01One}
(\d^k_{0,1}\vfi)(y_0,\ldots,y_{q})
=
\frac{2}{q+1}\,
\sum_{j=0}^{q}(-1)^{j}\,
m\left(
\vfi(y_0,\ldots,\widehat{y_j},\ldots,y_{q}),\,y_j
\right);
\end{equation}
if $p>0$, or if $q$ is even, then
\begin{equation}
\label{TheCobOp01}
\begin{array}{l}
(\d^k_{0,1}\vfi)(x_0,\ldots,x_{p-1};\,y_0,\ldots,y_{q})
=
\\[10pt]
\qquad
\qquad
\displaystyle
\frac{1}{q+1}\,
\sum_{j=0}^{q}(-1)^{p+j}\,
m\left(
\vfi(x_0,\ldots,x_{p-1};\,
y_0,\ldots,\widehat{y_j},\ldots,y_{q}),\,y_j
\right);
\end{array}
\end{equation}

(iii)
if $p>0$, then $\d^k_{-1,2}$ is given by
\begin{equation}
\label{TheCobOp-12}
\begin{array}{l}
(\d^k_{-1,2}\,\vfi)(x_0,\ldots,x_{p-2};\,y_0,\ldots,y_{q+1})
=
\\[10pt]
\qquad
\displaystyle
\frac{2}{(q+1)(q+2)}\,
\sum_{i<j}(-1)^{p+i+j}\,
\vfi(x_0,\ldots,x_{p-2},m(y_i,y_j);
\,y_0,\ldots,\widehat{y_i},\ldots,\widehat{y_j},\ldots,y_{q+1})
\end{array}
\end{equation}
is a coboundary operator, that is, $\d^{k+1}\circ\d^k=0$.
\end{cor}

 Note that the operator $\d_{1,0}$ is very close to the Hochschild differential,
 while $\d^k_{-1,2}$ is that of Chevalley-Eilenberg.
 
\begin{proof}
This is just the bracket $[m,.]_\mathfrak{al}$,
the above formul{\ae} is nothing but  (\ref{OurBr}) written explicitly.
We will give here the details in the most non-trivial case (ii).

Let first $p=0$ and $q$ be odd.
The (parity preserving) $(0,q)$-linear map $\vfi$ is with walues in $V_1$.
The formula of $\d_{01}$ is obtained from (\ref{GerProdGradTre}) by application of $\Alt$.
More precisely,
$$
j_{m}\,\vfi
\left(
y_0,\ldots,y_{q}
\right)=
m
\left(
\vfi(y_{0},\ldots,y_{q-1}),y_{q}
\right)+
m\left(
y_{0},
\vfi(y_{1},\ldots,y_{q})
\right).
$$
One therefore obtains:
$
\Alt
\left(j_{m}\,\vfi
\left(
y_0,\ldots,y_{q}
\right)\right)=
2\,\Alt\left({}m
\left(
\vfi(y_{0},\ldots,y_{q-1}),y_{q}
\right)\right)
$
in full accordance with (\ref{TheCobOp01One}).

If now $p>0$, then $\d_{01}$ is obtained by alternation of (\ref{GerProdGradBis}).
One has
$$
j_{m}\,\vfi
\left(
x_0,\ldots,x_{p-1};\,y_0,\ldots,y_{q}
\right)=
m
\left(
\vfi(x_{0},\ldots,x_{p-1};\,y_{0},\ldots,y_{q-1}),
y_{q}
\right).
$$
Applying the projector $\Alt$ one obtains (\ref{TheCobOp01}).
If finally $q$ is even, then $\vfi$ is with values in $V_0$ and the expression for $\d_{01}$
is again obtained by alternation of (\ref{GerProdGradBis}).
 \end{proof}

\begin{proposition}
The relation $\d^{k+1}\circ\d^k=0$ is equivalent to the following system:
\begin{equation}
\label{Delta2Eq}
\begin{array}{rcl}
\d^{k+1}_{1,0}\circ\d^{k}_{1,0} &=& 0 \\[6pt]
\d^{k+1}_{1,0}\circ\d^k_{0,1}+\d^{k+1}_{0,1}\circ\d^k_{1,0} &=& 0 \\[6pt]
\d^{k+1}_{0,1}\circ\d^k_{0,1}+
\d^{k+1}_{1,0}\circ\d^k_{-1,2} + \d^{k+1}_{-1,2}\circ\d^k_{1,0}&=& 0 \\[6pt]
\d^{k+1}_{0,1}\circ\d^k_{-1,2} + \d^{k+1}_{-1,2}\circ\d^k_{0,1}&=& 0 \\[6pt]
\d^{k+1}_{-1,2}\circ\d^k_{-1,2} &=& 0.
\end{array}
\end{equation}
\end{proposition}

\begin{proof}
The above equations obviously represent linearly independent
terms in
$\d^{k+1}\circ\d^k=0$.
\end{proof}

\subsection{Lie antialgebra and the zero-square condition}\label{AssSec}

Let us show that the notion of Lie antialgebra is related to graded Lie algebra
$\mathfrak{al}(V)$.

Given a Lie antialgebra $\fa$, in order to use the notation of Section \ref{GerAlt},
we denote by $V$ the ambient vector space, i.e.,
$V\cong\fa$.
Define the bilinear map $m:V\times{}V\to{}V$ as follows:
\begin{equation}
\label{TheOp}
\textstyle
m(x_1,x_2):=\frac{1}{2}\,x_1\cdot{}x_2\,,
\qquad
m(x,y):=x\cdot{}y\,,
\qquad
m(y_1,y_2):=y_1\cdot{}y_2\,,
\end{equation}
where $x_1,x_2\in{}V_0(=\fa_0)$ and $y_1,y_2\in{}V_1(=\fa_1)$
and where $\cdot$ stands for the product in $\fa$.

\begin{proposition}
\label{FormPro}
The operation $m$ satisfies the condition $[m,m]_\mathfrak{al}=0$
in the graded Lie algebra $\mathfrak{al}(V)$ if and only if $\fa$ is a
Lie antialgebra.
\end{proposition}

\begin{proof}
The condition $[m,m]_\mathfrak{al}=0$ reads:
$$
\begin{array}{rcl}
m\left(m(x_1,x_2),x_3\right)-
m\left(x_1,m(x_2,x_3)\right)&=&0,\\[6pt]
m\left(m(x_1,x_2),y\right)-
m\left(x_1,m(x_2,y)\right)&=&0,\\[6pt]
\textstyle
\half\,
m\left(m(x_1,y_1),y_2\right)-
\half\,m\left(m(x_1,y_2),y_1\right)
-m\left(x_1,m(y_1,y_2)\right)&=&0,\\[6pt]
m\left(m(y_1,y_2),y_3\right)+\hbox{cycle}
&=&0,
\end{array}
$$
cf. Example \ref{ImportEx}.
Using the definition (\ref{TheOp}), one immediately checks that the above condition
is indeed equivalent to the identities 
(\ref{AssCommT})--(\ref{Jack}).
\end{proof}

The structure of Lie antialgebra is therefore equivalent to the
zero-square condition in the graded Lie algebra $\mathfrak{al}(V)$.
It follows that the cohomology developed in Section \ref{TheDifSec}
is the cohomology of Lie antialgebras.
Let us give a detailed description of
cohomology of a Lie antialgebra with coefficients in an arbitrary module.

\subsection{Modules over Lie antialgebras}

So far, the cohomology we considered was with coefficients in the algebra itself.
However, this particular case actually contains the most general one.
The definition of a module is universal for all classes of algebras.

\begin{defi}
\label{MoDef}
{\rm
Let $(V,m)$ be an algebra (of an arbitrary type: associative, Lie, anti-Lie, etc.).
Assume that the vector space $V$ is a direct sum of two subspaces $V=V'\oplus{}W$.
If the space $V'$ is closed with respect to the bilinear map $m$,
in other words, 
$$
m:V'\times{}V'\to{}V',
$$
and, in addition, 
$$
m:V'\times{}W\to{}W
\qquad \hbox{and}\qquad
m:W\times{}W\to0,
$$
then the space $W$ is called a \textit{module} over the algebra $(V',m)$
}
\end{defi}

The space of multilinear maps from $V$ to $V$
contains the subspace, $C(V',W)$, of multilinear maps from $V'$ to $W$.
This subspace is obviously stable
under the differential $\d=\ad_m$.
The corresponding complex 
$$
\d:C^k(V',W)\to{}C^{k+1}(V',W)
$$
defines the cohomology of the algebra $(V',m)$ with coefficients
in the module $W$.

The notion of module over a Lie antialgebra fits into the general Definition \ref{MoDef}.
Given a Lie antialgebra $\fa$ and an $\fa$-module $\cB$, the space
$\fa\oplus\cB$ is equipped with a Lie antialgebra structure.
More precisely, for $a\in\fa$ and $b\in\cB$, one has
$$
(a,b)\cdot{}(a',b')=
\left(
a\cdot{}a',\,\r_ab'+(-1)^{\bar{a'}\bar{b}}\r_{a'}b
\right),
$$
where $\r:\fa\to\End(\cB)$ is the linear map
 that defines the $\fa$-action on $\cB$.
The Lie antialgebra structure on the space $\fa\oplus\cB$ is
called a \textit{semi-direct product} and is denoted by
$\fa\ltimes\cB$.

 \begin{proposition}\label{DualMod}
 Given an~$\fa$-module $\cB$, the dual space $\cB^*$
is naturally an~$\fa$-module the $\fa$-action being given by
\begin{equation}
\label{Diact}
\langle\r_a^*u,\,b\rangle:=(-1)^{\bar u\bar{a}}\,\langle u,\r_ab\rangle
\end{equation}
for all $u\in\cB^*$ and $b\in\cB$.
\end{proposition}
\begin{proof}
Straightforward.
\end{proof}

In particular, the space $\fa^*$ dual to a given Lie antialgebra $\fa$
 is an $\fa$-module, as the Lie antialgebra $\fa$ itself is obviously an $\fa$-module
with $\r_ab=a\cdot{}b$.
\subsection{Cohomology with coefficients in a module}\label{CohCoeff}

Given a Lie antialgebra $\fa$ and an $\fa$-module $\cB$,
we define the space, $C^{p,q}(\fa;\cB)$, of parity preserving skew-symmetric on $\fa_1$ maps
\begin{equation}
\label{CochEq}
\vfi:
\underbrace{
\left(
\fa_0\otimes\cdots\otimes\fa_0
\right)
}_{p}
\otimes
\underbrace{
\left(
\fa_1\wedge\cdots\wedge\fa_1
\right)
}_{q}
\to\cB.
\end{equation}
We also consider the following space:
$$
C^k(\fa;\cB)=
\bigoplus_{q+p=k}C^{q,p}(\fa;\cB)
$$
that we call the space of $k$-cochains.
One obviously has $C^k(\fa;\cB)\subset{}\mathfrak{al}(\fa\ltimes\cB)$.

The coboundary map $\d$ to the space of
cochains $C^{p,q}(\fa;\cB)$ is defined by formul{ae} 
(\ref{TheCobOp10})-(\ref{TheCobOp-12}) of Theorem \ref{CombThm},
where $m$ is as in (\ref{TheOp}).
We define the cohomology of a Lie antialgebra $\fa$ with
coefficients in an $\fa$-module $\cB$ in a usual way:
$$
H^k(\fa;\cB)=\ker(\d^k)/\mathrm{im}(\d^{k-1}),
$$
where $k\geq1$.
The space $\ker(\d^k)$ is called the space of $k$-cocycles
and the space $\mathrm{im}(\d^{k-1})$ is called the space of $k$-coboundaries.

\subsection{The case of trivial coefficients}

In the case, where $\cB=\bbR$ (or $\bbC$) is a
trivial module, the coboundary map (\ref{TheCobOpTot}) becomes simpler.
The operator $\d_{0,1}$ vanishes identically, while the system
(\ref{Delta2Eq}) reads:
$$
{\d_{1,0}}^2=0,
\qquad
{\d_{-1,2}}^2=0,
\qquad
\left[
\d_{1,0},\d_{-1,2}
\right]=0.
$$
One therefore obtains a structure of bicomplex
with two commuting differentials, $\d_{1,0}$ and $\d_{-1,2}$.
We denote by $H^k(\fa)$ the $k$-th cohomology space with trivial coefficients of
a Lie antialgebra $\fa$.

\begin{exe}
{\rm
Consider the Kaplansky superalgebra $K_3$.
In this case, the cohomology with trivial coefficients are
quite easy to calculate. This cohomology is trivial:
$H^k(K_3)=0$, for all $k>0$.
We omit the explicit computation.
}
\end{exe}

\section{Algebraic interpretation of lower degree cohomology}

Cohomology spaces of lower degree have algebraic
meaning quite similar to that in the usual case of Lie algebras.
In this section, we use the notation $a\cdot{}b$ for the action $\rho_ab$ of $a\in\fa$
on an element $b\in\cB$, thinking of $\cB$ as an ideal in $\fa\ltimes\cB$.

\subsection{The first cohomology space $H^1(\fa;\cB)$}\label{FirstSpace}

An even derivation of $\fa$ with values in the $\fa$-module $\cB$ is a
parity preserving linear map $c:\fa\to\cB$ such that
$$
c(a\cdot{}a')=c(a)\cdot{}a'+a\cdot{}c(a')
$$
where $\cdot$ stays both for the product in $\fa$ and for the $\fa$-action on $\cB$.

\begin{proposition}
\label{DerProp}
The first cohomology space $H^1(\fa;\cB)$ is the space of even derivations
of~$\fa$ with values in $\cB$.
\end{proposition}

\begin{proof}
A parity preserving linear map $c:\fa\to\cB$ is a sum
$c=c_{00}+c_{11}$, where $c_{00}:\fa_0\to\cB_0$ and $c_{11}:\fa_1\to\cB_1$.
The condition $\d\,c=0$ reads:
$$
\d_{10}\,c_{00}=0,
\qquad
\d_{01}\,c_{00}+\d_{10}\,c_{11}=0,
\qquad
\d_{-12}\,c_{00}+\d_{01}\,c_{11}=0.
$$
Let us show that these equations are equivalent to $c(a\cdot{}a')=c(a)\cdot{}a'+a\cdot{}c(a')$.

Since $\d_{10}\,c_{00}$ is nothing but the Hochshild coboundary
and $c_{00}$ is even, the first condition, $\d_{10}\,c_{00}(x_0,x_1)=0$, reads explicitly:
$-x_0\cdot{}c_{00}(x_1)-c_{00}(x_0)\cdot{}x_1+c_{00}(x_0\cdot{}x_1)=0$.

Furthermore, one obtains from (\ref{TheCobOp10}) and (\ref{TheCobOp01}):
$$
\d_{01}\,c_{00}(x,y)=-c_{00}(x)\cdot{}y,
\qquad
\d_{10}\,c_{11}(x,y)=-x\cdot{}c_{11}(y)+c_{11}(x\cdot{}y).
$$
The second equation is thus equivalent to 
$c_{00}(x)\cdot{}y+x\cdot{}c_{11}(y)-c_{11}(x\cdot{}y)=0$
and similar for the other two equations.

The space of coboundaries vanishes.
\end{proof}

\subsection{Second cohomology space $H^2(\fa,\cB)$ and abelian extensions}

An exact sequence of homomorphisms of Lie antialgebras
\begin{equation}
\label{SubmodDiagramm}
\begin{CD}
0 @> >>\cB @> >> \widetilde\fa @> >> \fa @> >>0,
\end{CD}
\end{equation}
where $\cB$ is a trivial algebra,
is called an abelian extension of the Lie antialgebra $\fa$
with coefficients in $\cB$.

As a vector space, $\widetilde\fa=\fa\oplus\cB$ and
the subspace $\cB$ is obviously an $\fa$-module,
the Lie antialgebra structure on $\widetilde\fa$ being given by
\begin{equation}
\label{CocEq}
(a,b)\cdot{}(a',b')=
\left(
a\cdot{}a',\,a\cdot{}b'+(-1)^{\bar{a'}\bar{b}}\,a'\cdot{}b+c(a,a')
\right),
\end{equation}
where $a,a'\in\fa$ and $b,b'\in\cB$ and where
$c:\fa\times\fa\to\cB$ is a bilinear map preserving the parity.

An extension~(\ref{SubmodDiagramm}) is called trivial if
the Lie antialgebra $\widetilde\fa$ is isomorphic to the semi-direct sum
$\fa\ltimes{}\cB$, i.e., to the extension (\ref{CocEq}) with $c=0$,
and the isomorphism is of the form
\begin{equation}
\label{IsomTriv}
(a,b)\mapsto{}(a,b+L(a)),
\end{equation}
where $L:\fa\to\cB$ is a Linear map.
If such an isomorphism does not exist, then the extension is called
\textit{non-trivial}.

The following statement shows that second cohomology space 
$H^2(\fa,\cB)$ classifies non-trivial extensions
of $\fa$ with coefficients in $\cB$.

\begin{proposition}
\label{CocProp}
(i)
The formula (\ref{CocEq}) defines a Lie antialgebra
structure if and only if the map $c$ is a 2-cocycle.

(ii) Two extensions $\widetilde\fa$ and ${\widetilde\fa}^\prime$
of the same Lie antialgebra $\fa$ by a module $\cB$ are
isomorphic
if and only if the corresponding cocycles $c_1$ and $c_2$
define the same cohomology class in  $H^2(\fa,\cB)$.
\end{proposition}

\begin{proof}
Part (i).
The axioms (\ref{AssCommT}-\ref{Jack}) applied to formula (\ref{CocEq}) read
\begin{equation}
\label{CocEqExp}
\begin{array}{rcl}
x_1\cdot{}c(x_2,x_3)
-c(x_1\cdot{}x_2,\,x_3)
+c(x_1,\,x_2\cdot{}x_3)
-c(x_1,x_2)\cdot{}x_3&=&0,\\[10pt]
c(x_1,x_2)\cdot{}y
+c(x_1\cdot{}x_2,\,y)
-2\,x_1\cdot{}c(x_2,y)
-2\,c(x_1,\,x_2\cdot{}y)&=&0,\\[10pt]
c(x,\,y_1\cdot{}y_2)
+x\cdot{}c(y_1,y_2)
-c(x,y_1)\cdot{}y_2&&\\[6pt]
-c(x\cdot{}\,y_1,\,y_2)
-c(y_1,\,x\cdot{}y_2)
-y_1\cdot{}c(x,y_2)
&=&0,\\[10pt]
y_1\cdot{}c(y_2,y_3)+c(y_1,\,y_2\cdot{}y_3)+(\hbox{cycle})&=&0.
\end{array}
\end{equation}
Substituting $m$ defined by (\ref{TheOp}),
one obtains precisely the condition $\d\,c=0$, where $\d$ is as in Theorem \ref{CombThm}.

Part (ii).
Assume that an extension (\ref{SubmodDiagramm}) is trivial
and there exists an isomorphism (\ref{IsomTriv}).
One then readily checks that this is equivalent to $c=\d\,L$.
\end{proof}

An extension (\ref{SubmodDiagramm}) of a Lie antialgebra is called a  central extension,
if the module $\cB$ is trivial, i.e., if the $\fa$-action on $\cB$ is identically zero.
In this case, $\cB$ belongs to the center of $\widetilde\fa$.
Central extensions with $\cB=\bbK$ are classified
by the second cohomology with trivial coefficients that we denote by
$H^2(\fa)$.

\section{Two remarkable cocycles}

In this section, we give examples of two non-trivial cohomology
classes of the infinite-dimensional Lie antialgebras
$\cAK(1)$ and $\cM^1$, see Examples \ref{ExMain} and \ref{ExMainBis}.
These cohomology classes are
analogues of the famous Gelfand-Fuchs and Godbillon-Vey
classes.

\subsection{The Gelfand-Fuchs cocycle}

Recall that the classical Gelfand-Fuchs cocycle is
a 2-cocycle with trivial coefficients on the Lie algebra of vector fields on the circle.
This cocycle defines the unique central extension of the Lie algebra of vector fields called
the Virasoro algebra.
This algebra plays an important r\^ole in geometry and
mathematical physics, see \cite{Fuk,GR} and references therein.
In the $\bbZ_2$-graded case, the graded version of the Gelfand-Fuchs cocycle
defined the Neveu-Schwarz and Ramond conformal superalgebras.

In this section, we recall the definition of the
Gelfand-Fuchs cocycle on the conformal Lie superalgebra $\cK(1)$
and rewrite it in a form of a 1-cocycle with coefficients in the dual space~$\cK(1)^*$.
This ``dualized'' version of the Gelfand-Fuchs cocycle will be of a particular interest
for our purpose.

The conformal Lie superalgebra, $\cK(1)$, is spanned by the basis
$$
\textstyle
\left\{
\ell_n,\;n\in\bbZ;
\qquad
\xi_i,\;i\in\bbZ+\half
\right\}
$$
with the following commutation relations
\begin{equation}
\label{CAlgRel}
\begin{array}{rcl}
\left[
\ell_n,\ell_m
\right] &=&
\left(m-n\right)\ell_{n+m},\\[8pt]
\left[
\ell_n,\xi_i
\right] &=&
\left(i-\frac{n}{2}\right)\xi_{n+i},\\[8pt]
\left[
\xi_i,\xi_j
\right] &=&
2\,\ell_{i+j}.
\end{array}
\end{equation}
The Lie subalgebra generated by $\ell_i$ is the Lie algebra of
(polynomial) vector fields on $S^1$.

The second cohomology space of $\cK(1)$ with trivial coefficients, $H^2(\cK(1))$, is
one-dimensional and is generated by the 2-cocycle
\begin{equation}
\label{GFClass}
\begin{array}{rcl}
c_{GF}(\ell_n,\ell_m)&=&
\left(
n^3-n
\right)
\d_{n+m,0}\\[5pt]
c_{GF}(\xi_i,\xi_j)&=&
\left(
-4\,i^2+1
\right)
\d_{i+j,0}\\[5pt]
c_{GF}(\ell_n,\xi_i)&=&
0,
\end{array}
\end{equation}
that we call the Gelfand-Fuchs cocycle.
It defines a (unique) central extension of $\cK(1)$.
The even part of $c_{GF}$ defines the Virasoro algebra.

\begin{rmk}
{\rm
The conformal Lie superalgebra $\cK(1)$ contains a subalgebra spanned by the elements
$$
\left\{
\ell_{-1},\,\ell_0,\,\ell_1;\,\xi_{-\half},\,\xi_\half
\right\},
$$
isomorphic to  the Lie superalgebra $\mathrm{osp}(1|2)$.
The cocycle (\ref{GFClass}) can be characterized as
\textit{the unique} 2-cocycle on $\cK(1)$ vanishing on $\mathrm{osp}(1|2)$.
}
\end{rmk}

\subsection{The dual Gelfand-Fuchs cocycle}\label{Dus}

It is a general fact that a 2-cocycle on a Lie superalgebra with trivial coefficients
corresponds to a 1-cocycle with coefficients in the dual space.

\begin{defi}
{\rm
Given a Lie algebra $\fg$ and a 2-cocycle $c:\fg\wedge\fg\to\bbK$,
the formula
$$
\langle{}C(X),\,Y\rangle:=c(X,Y),
$$
for all $X,Y\in\fg$, defines a
1-cocycle $C:\fg\to\fg^*$ with values in the dual space $\fg^*$.
}
\end{defi}

The ``dual Gelfand-Fuchs cocycle'' then reads:
\begin{equation}
\label{VraiGF}
\textstyle
C_{GF}(\ell_n)=(n^3-n)\,\ell^*_{-n},
\qquad
C_{GF}(\xi_i)=(-4\,i^2+1)\,\xi^*_{-i},
\end{equation}
where $\{\ell^*_{n},\,\xi^*_{i}\}$ is the dual basis, i.e.,
$$
\langle\ell^*_{n},\ell_m\rangle=\d_{n,m},
\qquad
\langle\xi^*_i,\xi_j\rangle=\d_{i,j},
\qquad
\langle\ell^*_{n},\xi_j\rangle=\langle\xi^*_i,\ell_m\rangle=0.
$$
This is a non-trivial 1-cocycle on $\cK(1)$ with values in $\cK(1)^*$.

\subsection{A non-trivial cocycle on $\cAK(1)$}

The conformal Lie antialgebra $\cAK(1)$ has no non-trivial
central extension, see \cite{Ovs},
therefore, there is no analog of the classical Gelfand-Fuchs cocycle (\ref{GFClass}).
However, there exists an analog of the dual cocycle (\ref{VraiGF}).
We denote by $\{\e_n^*,a_i^*\}$ the of $\cAK(1)^*$ dual to $\{\e_n,a_i\}$.

\begin{theorem}
The linear map $\g:\cAK(1)\to\cAK(1)^*$ given by
\begin{equation}
\label{FauxGF}
\textstyle
\g(\e_n)=-n\,\e_{-n}^*,
\qquad
\g(a_i)=\left(i^2-\frac{1}{4}\right)a_{-i}^*,
\end{equation}
is a non-trivial 1-cocycle on $\cAK(1)$.
\end{theorem}

\begin{proof}
Let us check that the map (\ref{FauxGF}) is, indeed, a 1-cocycle.
The action of $\cAK(1)$ on $\cAK(1)^*$ can be easily calculated according to the formula (\ref{Diact}).
The result is as follows:
$$
\begin{array}{rcl}
\e_n\cdot{}\e_m^*&=&\e_{m-n}^*,\\[6pt]
\e_n\cdot{}a_i^*&=&\half\,a_{i-n}^*\\[6pt]
a_i\cdot{}\e_m^*&=&(\frac{m}{2}-i)\,a_{m-i}^*,\\[6pt]
a_i\cdot{}a_j^*&=&-\half\,\e_{j-i}^*.
\end{array}
$$
Consider the following Ansatz:
$$
\g(\e_n)=t(n)\,\e^*_{-n},
\qquad
\g(a_i)=s(i)\,a^*_{-i}.
$$
The 1-cocycle condition then leads to: 
$\g(\e_n\cdot\e_m)=\e_n\cdot\g(e_m)+\e_m\cdot\g(e_n)$, so that
$$
t(n+m)=t(n)+t(m).
$$
The next two conditions are:
$\half\g(\e_n\cdot{}a_i)=\e_n\cdot\g(a_i)+\g(\e_n)\cdot{}a_i$ and
$\g(a_i\cdot{}a_j)=a_i\cdot\g(a_j)-a_j\cdot\g(a_i)$.
They give the same equation:
$$
s(i)-s(j)=i^2-j^2.
$$
The map (\ref{FauxGF}) obviously satisfies both equations, so that this is, indeed, a 1-cocycle.

We have already seen in Section \ref{FirstSpace} that every coboundary
vanishes on the even part of a Lie antialgebra.
It follows that the 1-cocycle (\ref{FauxGF}) is non-trivial.
\end{proof}

We conjecture that the space $H^1(\cAK(1);\,\cAK(1)^*)$ is one-dimensional
and, thus generated by the 1-cocycle (\ref{FauxGF}).
We think that the cocycle (\ref{FauxGF}) is characterized by the property that
it vanishes on the subalgebra $K_3\subset\cAK(1)$ spanned by $\e_0,a_{-\half},a_\half$.

\begin{rmk}
{\rm
The 1-cocycle (\ref{FauxGF}) is skew-symmetric on the even part $\cAK(1)_0$
and symmetric on the odd part $\cAK(1)_1$.
This is the reason why it cannot be understood as a 2-cocycle on $\cAK(1)$
with trivial coefficients.
This phenomenon is quite general for the Lie antialgebras:
approximately a half of the statements that hold in the Lie algebra setting
remains true.
}
\end{rmk}

\subsection{A cocycle on $\cM^1$ and the dual Godbillon-Vey cocycle}

The Godbillon-Vey cocycle is a 3-cocycle on the Lie algebra of polynomial
vector fields on the line.
The simplest way to define this cocycle is as follows.
Consider the Lie algebra $W_1$ with basis $\{\ell_{n},\;n\geq-1\}$
and the relations $[\ell_n,\ell_m]=(m-n)\,\ell_{n+m}$.
Then
$$
c_{GV}=\ell^*_{-1}\wedge\ell^*_0\wedge\ell^*_1
$$
is a non-trivial 3-cocycle on $W_1$ (with trivial coefficients).
Similarly to Section \ref{Dus}, one can ``dualize'' the above 3-cocycle and
obtain a 2-cocycle with coefficients in $W^*_1$.

Let us consider the Lie antialgebra $\cM^1$, see Example \ref{ExMainBis}.

\begin{theorem}
The bilinear map $\eta:\cM^1\otimes\cM^1\to(\cM^1)^*$ given by
\begin{equation}
\label{FauxGV}
\textstyle
\eta=a_{-\half}^*\wedge{}a_{\half}^*\otimes\e^*_0,
\end{equation}
is a non-trivial 2-cocycle on $\cM^1$.
\end{theorem}

\begin{proof}
First, one can easily check that
$$
\textstyle
\eta'=\l\,a_{-\half}^*\wedge{}a_{\half}^*\otimes\e^*_0+
\mu\left(
\e_0^*\wedge\e_0^*\otimes\e^*_0-
\half\,\e_0^*\wedge a_{-\half}^*\otimes{}a^*_\half+
\half\,\e_0^*\wedge{}a_{\half}^*\otimes{}a^*_{-\half}
\right)
$$
is a 2-cocycle for all $\l,\mu$.

Second, for any coboundary $\d\zeta$, where $\zeta:\cM^1\to(\cM^1)^*$ is a linear map,
one has
$$
(\d\zeta)(\e_0,\e_0)=2\,(\d\zeta)(a_{-\half},a_{\half}),
$$
so that $\l=\half\mu$ if and only if $\eta'$ is a coboundary.
\end{proof}

We conjecture that the space $H^2(\cM^1;(\cM^1)^*)$ is one-dimensional.
We also hope that the cocycle (\ref{FauxGV}) has a topological meaning
and can be associated to a characteristic class, similarly to the classical
Godbillon-Vey cocycle, see \cite{Fuk}.

\medskip

\noindent \textbf{Acknowledgments}.
We are grateful to
F. Chapoton, S. Leidwanger, J.-L. Loday, S. Morier-Genoud
for enlightening discussions
and to the referee for helpful comments.
Special thanks are due to Marie Kreusch for a careful reading of
a preliminary version of this paper.

\vskip 1cm


Pierre Lecomte, {\hfill Valentin Ovsienko,}

D\'epartement de math\'ematiques, {\hfill CNRS,}

Universit\'e de Li\`ege,   {\hfill Institut Camille Jordan,}

Grande Traverse 12 (B37), {\hfill Universit\'e Lyon 1,}

B-4000 Li\`ge, {\hfill Villeurbanne Cedex, F-69622,}

Belgique; {\hfill France;}

plecomte@ulg.ac.be {\hfill ovsienko@math.univ-lyon1.fr}

\end{document}